\documentclass[a4paper,12pt]{article} 
\usepackage{amsmath} 
\usepackage{amsopn} 
\usepackage{amssymb} 
\usepackage[mathscr]{eucal} 
\usepackage{theorem} 
\usepackage{enumerate} 

\theorembodyfont{\rmfamily} 
  \newtheorem{rem}{Remark}

\theorembodyfont{\itshape} 
\theoremstyle{plain}
  \newtheorem{thm}{Theorem}[section] 
  \newtheorem{pro}[thm]{Proposition} 
   
  \newtheorem{cor}[thm]{Corollary} 

\renewcommand{\theequation}%
           {\thesection.\arabic{equation}}

\begin{document} 

\begin{center} 
{\LARGE The $SO(3,1)$-orbits in the light cone of} 

\vspace{1mm} 

{\LARGE the 2-fold exterior power of the Minkowski 4-space} 

\vspace{6mm} 

{\Large Naoya {\sc Ando}} 
\end{center} 

\vspace{3mm} 

\begin{quote} 
{\footnotesize \it Abstract} \ 
{\footnotesize Two special neutral hypersurfaces $\mathcal{L}_{\pm}$ 
in the light cone $L(\bigwedge^2 E^4_1 )$ studied in \cite{ando}, \cite{ando'} 
are $SO(3,1)$-orbits. 
In this paper, we see that each $SO(3,1)$-orbit in $L(\bigwedge^2 E^4_1 )$ 
is either a neutral hypersurface homothetic to one of $\mathcal{L}_{\pm}$ 
in $L(\bigwedge^2 E^4_1 )$ 
or a hypersurface with a two-dimensional involutive distribution 
where the induced metric is degenerate. 
The difference between these hypersurfaces can be understood 
in terms of the stabilizer and 
the $r$-slice of $L(\bigwedge^2 E^4_1 )$ for $r>0$.} 
\end{quote} 

\vspace{3mm} 

\section{Introduction}\label{sect:intro} 

\setcounter{equation}{0} 

The purpose of this paper is to study 
the $SO(3,1)$-orbits in the light cone $L(\bigwedge^2 E^4_1 )$ 
of $\bigwedge^2 E^4_1$. 

The $SO(3,1)$-action on the Minkowski $4$-space $E^4_1$ yields 
a natural $SO(3,1)$-action on the $2$-fold exterior power $\bigwedge^2 E^4_1$ 
of $E^4_1$. 
In particular, $SO(3,1)$ acts on the light cone $L(\bigwedge^2 E^4_1 )$ 
of $\bigwedge^2 E^4_1$. 
Two special $SO(3,1)$-orbits in $L(\bigwedge^2 E^4_1 )$ appear in 
studies of space-like or time-like surfaces in $E^4_1$ 
with zero mean curvature vector (\cite{ando}). 
Let $M$ be a Riemann surface 
and $F:M\longrightarrow E^4_1$ a space-like and conformal immersion of $M$ 
into $E^4_1$ with zero mean curvature vector. 
Then the lifts $\Omega_{F, \pm}$ of $F$ are maps from $M$ 
into neutral hypersurfaces $\mathcal{L}_{\pm}$ of $L(\bigwedge^2 E^4_1 )$ 
respectively, 
and they are holomorphic with respect to 
parallel almost complex structures of $\mathcal{L}_{\pm}$ 
(\cite{ando}, \cite{ando'}). 
If we suppose that $M$ is a Lorentz surface 
and that $F:M\longrightarrow E^4_1$ is a time-like and conformal immersion 
of $M$ into $E^4_1$ with zero mean curvature vector, 
then we have $\Omega_{F, \pm} :M\longrightarrow \mathcal{L}_{\pm}$ and 
they are holomorphic with respect to parallel almost paracomplex structures 
of $\mathcal{L}_{\pm}$ (\cite{ando}, \cite{ando'}). 
The hypersurfaces $\mathcal{L}_{\pm}$ are $SO(3,1)$-orbits 
in $L(\bigwedge^2 E^4_1 )$. 

In this paper, we will see that each $SO(3,1)$-orbit $\mathcal{L}$ 
in $L(\bigwedge^2 E^4_1 )$ is 
either a neutral hypersurface 
which is homothetic to one of $\mathcal{L}_{\pm}$ 
in $L(\bigwedge^2 E^4_1 )$ (Proposition~\ref{pro:ts}, Theorem~\ref{thm:gcase}) 
or a hypersurface with a two-dimensional involutive distribution 
where the induced metric is degenerate (Theorem~\ref{thm:Dl}). 
In order to understand the difference between these hypersurfaces, 
we will study the stabilizer of $SO(3, 1)$ at a point of $\mathcal{L}$ 
and the invariant subspaces of the tangent space by the stabilizer 
(Proposition~\ref{pro:invsubsp}, Theorem~\ref{thm:stab}). 
In addition, we will define the $r$-slice $L_r (\bigwedge^2 E^4_1 )$ 
and study the intersection of the $r$-slice and each $SO(3,1)$-orbit 
in Section~\ref{sect:rslice}. 
The $r$-slice is a subset of $L(\bigwedge^2 E^4_1 )$ 
defined for each $r>0$ and 
each equivalence class in the set $\mathcal{B} (E^4_1 )$ 
of ordered pseudo-orthonormal bases of $E^4_1$ giving the orientation 
with respect to an equivalence relation given in Section~\ref{sect:lc}. 
For each neutral $SO(3,1)$-orbit $\mathcal{L}$ 
and each equivalence class in $\mathcal{B} (E^4_1 )$, 
we will see that there exists a positive number $r_0 >0$ such that 
the intersection $\mathcal{L} \cap L_r (\bigwedge^2 E^4_1 )$ is not empty 
if and only if $r\geq r_0$ and that 
if $\mathcal{L} \cap L_r (\bigwedge^2 E^4_1 )\not= \emptyset$, 
then $\mathcal{L} \cap L_r (\bigwedge^2 E^4_1 )$ is 
diffeomorphic to $\mbox{\boldmath{$R$}}P^3$ or $S^2$ 
according to $r>r_0$ or $r=r_0$ (Theorem~\ref{thm:slice1}). 
We will see that if an $SO(3,1)$-orbit $\mathcal{L}$ has 
a two-dimensional involutive distribution 
where the induced metric is degenerate, 
then $\mathcal{L} \cap L_r (\bigwedge^2 E^4_1 )$ is not empty, 
and diffeomorphic to $\mbox{\boldmath{$R$}}P^3$ for any $r>0$ 
(Theorem~\ref{thm:slice2}). 

\begin{rem} 
Holomorphicity of the lifts of space-like or time-like surfaces in $E^4_1$ 
with zero mean curvature vector is an analogue of 
holomorphicity of the Gauss maps of minimal surfaces in 
the Euclidean $4$-space $E^4$ (\cite[pp.\,16--22]{HO}) 
and holomorphicity of the Gauss maps of 
space-like or time-like surfaces with zero mean curvature vector 
in $E^4_2$ (\cite{ando2}). 
\end{rem} 

\begin{rem} 
In Proposition 1 of \cite{ando}, 
it was asserted that $\mathcal{L}_{\pm}$ are flat. 
However, this is not true. 
In \cite{ando'}, the curvature tensors of $\mathcal{L}_{\pm}$ are 
explicitly represented and in particular, they do not vanish. 
\end{rem} 

\section{The light cone of 
\mbox{\boldmath{$\bigwedge^2 E^4_1$}}}\label{sect:lc} 

\setcounter{equation}{0} 

Let $h$ be the metric of the Minkowski $4$-space $E^4_1$: 
$$h(x, y):=x^1 y^1 +x^2 y^2 +x^3 y^3 -x^4 y^4$$ 
for $x={}^t (x^1 \ x^2 \ x^3 \ x^4)$, 
    $y={}^t (y^1 \ y^2 \ y^3 \ y^4)\in E^4_1$. 
Let $\hat{h}$ be the metric of the $2$-fold exterior power $\bigwedge^2 E^4_1$ 
of $E^4_1$ induced by $h$: 
$$\hat{h} (x\wedge y, u\wedge v):=h(x, u)h(y, v)-h(x, v)h(y, u) \quad 
(x, y, u, v\in E^4_1 ).$$ 
Let $(e_1 , e_2 , e_3 , e_4 )$ be an element of $\mathcal{B} (E^4_1 )$, and 
suppose that $e_1$, $e_2$, $e_3$ are space-like and that $e_4$ is time-like. 
Then $\omega_{ij} :=e_i \wedge e_j$ ($1\leq i<j\leq 4$) form 
a pseudo-orthonormal basis of $\bigwedge^2 E^4_1$ with respect to $\hat{h}$ 
satisfying 
\begin{itemize} 
\item[{\rm (i)}]{$\omega_{12}$, $\omega_{13}$, $\omega_{23}$ are 
space-like;} 
\item[{\rm (ii)}]{$\omega_{14}$, $\omega_{24}$, $\omega_{34}$ are 
time-like.} 
\end{itemize}  
In particular, 
we see that $\bigwedge^2 E^4_1$ is of dimension $6$ 
and that $\hat{h}$ has signature $(3,3)$. 
For $\Omega =\sum_{1\leq i<j\leq 4} c_{ij} \omega_{ij}$, 
we set 
\begin{equation*} 
A(\Omega )
:=\sum_{1\leq i<j\leq 3} c^2_{ij} , \quad 
B(\Omega )
:=\sum_{1\leq i\leq 3} c^2_{i4} . 
\end{equation*} 
The light cone $L(\bigwedge^2 E^4_1 )$ of $\bigwedge^2 E^4_1$ is given by 
$$L(\textstyle\bigwedge^2 E^4_1 ) 
:=\left\{ \left. 
  \Omega \in \textstyle\bigwedge^2 E^4_1 
  \ \right| \ 
  A(\Omega )=B(\Omega )\not= 0\right\} .$$ 
We see that $L(\bigwedge^2 E^4_1 )$ does not depend on the choice 
of $(e_1 , e_2 , e_3 , e_4 )\in \mathcal{B} (E^4_1 )$. We set 
$$L_1 (\textstyle\bigwedge^2 E^4_1 ) 
:=\left\{ \left. 
  \Omega \in L(\textstyle\bigwedge^2 E^4_1 ) 
  \ \right| \ 
  A(\Omega )=B(\Omega )=1\right\} .$$ 
Then $L_1 (\bigwedge^2 E^4_1 )$ depends on the choice 
of $(e_1 , e_2 , e_3 , e_4 )\in \mathcal{B} (E^4_1 )$. 
For elements $e =(e_1  , e_2  , e_3  , e_4 )$, 
             $e'=(e'_1 , e'_2 , e'_3 , e'_4 )$ 
of $\mathcal{B} (E^4_1 )$, we write $e\sim e'$ if 
$$(e'_1 \ e'_2 \ e'_3 )=\varepsilon (e_1  \ e_2  \ e_3 )U, \quad 
   e'_4                =\varepsilon    e_4$$ 
for an element $U$ of $SO(3)$ and $\varepsilon \in \{ +, -\}$. 
Then $\sim$ is an equivalence relation in $\mathcal{B} (E^4_1 )$. 
Let $\Omega =\sum_{i<j} c_{ij} \omega_{ij}$ be an element 
of $L_1 (\bigwedge^2 E^4_1 )$ 
for $e =(e_1 , e_2 , e_3 , e_4 )\in \mathcal{B} (E^4_1 )$. 
We set 
\begin{equation*} 
\begin{array}{lclcl} 
 a_1 :=c_{23} , & \ & a_2 :=-c_{13} , & \ & a_3 :=c_{12} , \\ 
 b_1 :=c_{14} , & \ & b_2 := c_{24} , & \ & b_3 :=c_{34} 
  \end{array} 
\end{equation*} 
and $\mbox{\boldmath{$a$}} :={}^t (a_1 \ a_2 \ a_3 )$, 
    $\mbox{\boldmath{$b$}} :={}^t (b_1 \ b_2 \ b_3 )$. 
Then the lengths of $\mbox{\boldmath{$a$}}$, $\mbox{\boldmath{$b$}}$ are 
equal to one and $\Omega$ is represented as 
\begin{equation} 
\Omega =(\omega_{23} \ \omega_{31} \ \omega_{12} )\mbox{\boldmath{$a$}} 
       +(\omega_{14} \ \omega_{24} \ \omega_{34} )\mbox{\boldmath{$b$}} . 
\label{Omega} 
\end{equation} 
If $e\sim e'$, then we see by \eqref{Omega} that $\Omega$ is represented as 
\begin{equation} 
\Omega =(\omega'_{23} \ \omega'_{31} \ \omega'_{12} )\,{}^tU 
         \mbox{\boldmath{$a$}} 
       +(\omega'_{14} \ \omega'_{24} \ \omega'_{34} )\,{}^tU 
         \mbox{\boldmath{$b$}} , 
\label{Omega'} 
\end{equation} 
where $\omega'_{ij} :=e'_i \wedge e'_j$. 
Therefore $L_1 (\bigwedge^2 E^4_1 )$ is determined by 
each equivalence class with respect to $\sim$ 
and therefore we denote it by $L_1 (\bigwedge^2 E^4_1 , [e])$ 
for the equivalence class $[e]$ containing $e\in \mathcal{B} (E^4_1 )$. 

\begin{pro}\label{pro:L1} 
For each element $\Omega$ of $L_1 (\bigwedge^2 E^4_1 , [e])$, 
there exists an element $e'$ of $[e]$ satisfying 
\begin{equation*} 
\Omega =(\cos \phi )\omega'_{12} +(\sin \phi )\omega'_{23} +\omega'_{34} 
\end{equation*} 
for a number $\phi \in [0, \pi ]$. 
\end{pro} 

\vspace{3mm} 

\par\noindent 
\textit{Proof} \ 
We suppose that $\Omega$ is represented as in \eqref{Omega}. 
We represent $U$ as $U=(\mbox{\boldmath{$u$}}_1 \ 
                        \mbox{\boldmath{$u$}}_2 \ 
                        \mbox{\boldmath{$u$}}_3 )$. 
If we set $\mbox{\boldmath{$u$}}_3 =\mbox{\boldmath{$b$}}$ and 
if we suppose $\langle \mbox{\boldmath{$u$}}_2 , \mbox{\boldmath{$a$}} \rangle 
               =0$ for the standard inner product $\langle \ , \ \rangle$ 
of $\mbox{\boldmath{$R$}}^3$, 
then we see from \eqref{Omega'} that $\Omega$ is represented as 
\begin{equation*} 
 \Omega 
=\langle \mbox{\boldmath{$u$}}_1 , \mbox{\boldmath{$a$}} \rangle 
 \omega'_{23} 
+\langle \mbox{\boldmath{$u$}}_3 , \mbox{\boldmath{$a$}} \rangle 
 \omega'_{12} 
+\omega'_{34} . 
\end{equation*} 
Since $\langle \mbox{\boldmath{$u$}}_2 , \mbox{\boldmath{$a$}} \rangle =0$, 
we have $\langle \mbox{\boldmath{$u$}}_1 , \mbox{\boldmath{$a$}} \rangle^2 
        +\langle \mbox{\boldmath{$u$}}_3 , \mbox{\boldmath{$a$}} \rangle^2 
        =1$. 
Therefore there exists a number $\phi \in \mbox{\boldmath{$R$}}$ 
satisfying $\cos \phi =\langle \mbox{\boldmath{$u$}}_3 , 
                               \mbox{\boldmath{$a$}} \rangle$ 
and        $\sin \phi =\langle \mbox{\boldmath{$u$}}_1 , 
                               \mbox{\boldmath{$a$}} \rangle$. 
We can choose $\mbox{\boldmath{$u$}}_1$ 
satisfying $\langle \mbox{\boldmath{$u$}}_1 , \mbox{\boldmath{$a$}} \rangle 
            \geq 0$. 
Hence we obtain Proposition~\ref{pro:L1}. 
\hfill 
$\square$ 

\vspace{3mm} 

From Proposition~\ref{pro:L1}, we see that 
for each element $\Omega$ of $L(\bigwedge^2 E^4_1 )$, 
there exist an element $e=(e_1 , e_2 , e_3 , e_4 )$ 
of $\mathcal{B} (E^4_1 )$, $r>0$ and $\phi \in [0, \pi ]$ satisfying 
\begin{equation} 
    \Omega 
=r((\cos \phi )\omega_{12} +(\sin \phi )\omega_{23} +\omega_{34} ). 
\label{Omegar} 
\end{equation} 

\section{The \mbox{\boldmath{$SO(3,1)$}}-action 
on \mbox{\boldmath{$L(\bigwedge^2 E^4_1 )$}}}\label{sect:SO(3,1)action} 

\setcounter{equation}{0} 

For $P\in SO(3,1)$, 
we set $\Tilde{T}_P (\omega_{ij} ):=Pe_i \wedge Pe_j$. 
For an element $\Omega =\sum_{1\leq i<j\leq 4} c_{ij} \omega_{ij}$ 
of $\bigwedge^2 E^4_1$, we set 
\begin{equation*} 
  \Tilde{T}_{P} (\Omega ) 
:=\sum_{1\leq i<j\leq 4} c_{ij} \Tilde{T}_P (\omega_{ij} ). 
\end{equation*} 
Then we obtain a linear transformation $\Tilde{T}_{P}$ of $\bigwedge^2 E^4_1$, 
which does not depend on the choice 
of $(e_1 , e_2 , e_3 , e_4 )\in \mathcal{B} (E^4_1 )$. 
Since $\Tilde{T}_{P} \circ \Tilde{T}_{P'} =\Tilde{T}_{PP'}$, 
we obtain an $SO(3,1)$-action on $\bigwedge^2 E^4_1$. 
We see that $\Tilde{T}_{P}$ is an isometry of $\bigwedge^2 E^4_1$ 
with respect to the metric $\hat{h}$. 
In particular, if $\Omega \in L(\bigwedge^2 E^4_1 )$, 
then $\Tilde{T}_{P} (\Omega ) \in L(\bigwedge^2 E^4_1 )$, 
and therefore we obtain an $SO(3,1)$-action on $L(\bigwedge^2 E^4_1 )$. 
Suppose that $\Omega$ is represented as in \eqref{Omegar}. 
The $SO(3,1)$-orbit through $\Omega$ is given by 
\begin{equation} 
\mathcal{L} (\Omega )=\{ \Tilde{T}_{P} (\Omega ) \ | \ P\in SO(3,1)\} . 
\label{HOmega} 
\end{equation} 
Although $SO(3,1)$ has just two connected components, 
$\mathcal{L} (\Omega )$ is connected. 
In the following, we suppose $r=\sqrt{2}$: 
\begin{equation} 
 \Omega 
=\sqrt{2}\,((\cos \phi )\omega_{12} +(\sin \phi )\omega_{23} +\omega_{34} ). 
 \label{Omegasqrt2} 
\end{equation} 
As in \cite{ando}, we set 
\begin{equation*} 
\begin{split} 
& P_{1, 1} :=\left[ 
             \begin{array}{cccc} 
              \cos \theta & -\sin \theta & 0 & 0 \\ 
              \sin \theta &  \cos \theta & 0 & 0 \\ 
               0          &   0          & 1 & 0 \\ 
               0          &   0          & 0 & 1 
               \end{array} 
             \right] , \ 
  P_{1, 2} :=\left[ 
             \begin{array}{cccc} 
              1          &  0          &   0           &   0 \\ 
              0          &  1          &   0           &   0 \\ 
              0          &  0          & {\rm cosh}\,t & {\rm sinh}\,t \\ 
              0          &  0          & {\rm sinh}\,t & {\rm cosh}\,t 
               \end{array} 
             \right] , \\ 
& P_{2, 1} :=\left[ 
             \begin{array}{cccc} 
              \cos \theta & 0 & -\sin \theta & 0 \\ 
               0          & 1 &   0          & 0 \\ 
              \sin \theta & 0 &  \cos \theta & 0 \\ 
               0          & 0 &   0          & 1 
               \end{array} 
             \right] , \ 
  P_{2, 2} :=\left[ 
             \begin{array}{cccc} 
              1          &   0           &   0           &   0           \\ 
              0          & {\rm cosh}\,t &   0           & {\rm sinh}\,t \\ 
              0          &   0           &   1           &   0           \\ 
              0          & {\rm sinh}\,t &   0           & {\rm cosh}\,t 
               \end{array} 
             \right] , \\ 
& P_{3, 1} :=\left[ 
             \begin{array}{cccc} 
               1          &  0          &   0          & 0 \\ 
               0          & \cos \theta & -\sin \theta & 0 \\ 
               0          & \sin \theta &  \cos \theta & 0 \\ 
               0          &  0          &   0          & 1 
               \end{array} 
             \right] , \ 
  P_{3, 2} :=\left[ 
             \begin{array}{cccc} 
              {\rm cosh}\,t & 0          & 0             & {\rm sinh}\,t \\ 
                0           & 1          & 0             &   0           \\ 
                0           & 0          & 1             &   0           \\ 
              {\rm sinh}\,t & 0          & 0             & {\rm cosh}\,t 
               \end{array} 
             \right] 
\end{split} 
\end{equation*} 
for $\theta , t\in \mbox{\boldmath{$R$}}$. 
Then the connected component $SO_0 (3,1)$ of the unit element of $SO(3,1)$ is 
generated by $P_{k, l}$ ($k=1, 2, 3$, $l=1, 2$, 
$\theta , t\in \mbox{\boldmath{$R$}}$). 
We set 
\begin{equation*} 
\begin{split} 
E_{\pm , 1} & :=\dfrac{1}{\sqrt{2}} (\omega_{12} \pm \omega_{34} ), \\ 
E_{\pm , 2} & :=\dfrac{1}{\sqrt{2}} (\omega_{13} \pm \omega_{42} ), \\ 
E_{\pm , 3} & :=\dfrac{1}{\sqrt{2}} (\omega_{14} \pm \omega_{23} ). 
\end{split} 
\end{equation*} 
Then we have 
\begin{equation} 
\begin{split} 
&  [\Tilde{T}_{P_{1, 1}} (E_{\pm , 1} ) \ 
    \Tilde{T}_{P_{1, 1}} (E_{\pm , 2} ) \ 
    \Tilde{T}_{P_{1, 1}} (E_{\pm , 3} )] \\ 
& =[ E_{\pm , 1} \ E_{\pm , 2} \ E_{\pm , 3} ] 
    \left[ 
    \begin{array}{ccc} 
     1 &      0          &      0          \\ 
     0 &     \cos \theta & \mp \sin \theta \\ 
     0 & \pm \sin \theta &     \cos \theta 
      \end{array} 
    \right] , \\ 
&  [\Tilde{T}_{P_{1, 2}} (E_{\pm , 1} ) \ 
    \Tilde{T}_{P_{1, 2}} (E_{\pm , 2} ) \ 
    \Tilde{T}_{P_{1, 2}} (E_{\mp , 3} )] \\ 
& =[ E_{\pm , 1} \ E_{\pm , 2} \ E_{\mp , 3} ] 
    \left[ 
    \begin{array}{ccc} 
     1 &   0           &   0           \\ 
     0 & {\rm cosh}\,t & {\rm sinh}\,t \\ 
     0 & {\rm sinh}\,t & {\rm cosh}\,t 
      \end{array} 
    \right] , 
\end{split} 
\label{P11P12E} 
\end{equation} 
\begin{equation} 
\begin{split} 
&  [\Tilde{T}_{P_{2, 1}} (E_{\pm , 1} ) \ 
    \Tilde{T}_{P_{2, 1}} (E_{\pm , 2} ) \ 
    \Tilde{T}_{P_{2, 1}} (E_{\pm , 3} )] \\ 
& =[ E_{\pm , 1} \ E_{\pm , 2} \ E_{\pm , 3} ] 
    \left[ 
    \begin{array}{ccc} 
         \cos \theta & 0 & \pm \sin \theta \\ 
          0          & 1 &      0          \\ 
     \mp \sin \theta & 0 &     \cos \theta 
      \end{array} 
    \right] , \\ 
&  [\Tilde{T}_{P_{2, 2}} (E_{\pm , 1} ) \ 
    \Tilde{T}_{P_{2, 2}} (E_{\pm , 2} ) \ 
    \Tilde{T}_{P_{2, 2}} (E_{\mp , 3} )] \\ 
& =[ E_{\pm , 1} \ E_{\pm , 2} \ E_{\mp , 3} ] 
    \left[ 
    \begin{array}{ccc} 
     {\rm cosh}\,t & 0 & {\rm sinh}\,t \\ 
       0           & 1 &   0           \\ 
     {\rm sinh}\,t & 0 & {\rm cosh}\,t 
      \end{array} 
    \right] 
\end{split} 
\label{P21P22E} 
\end{equation} 
and 
\begin{equation} 
\begin{split} 
&  [\Tilde{T}_{P_{3, 1}} (E_{\pm , 1} ) \ 
    \Tilde{T}_{P_{3, 1}} (E_{\pm , 2} ) \ 
    \Tilde{T}_{P_{3, 1}} (E_{\pm , 3} )] \\ 
& =[ E_{\pm , 1} \ E_{\pm , 2} \ E_{\pm , 3} ] 
    \left[ 
    \begin{array}{ccc} 
     \cos \theta & -\sin \theta & 0 \\ 
     \sin \theta &  \cos \theta & 0 \\ 
      0          &   0          & 1 
      \end{array} 
    \right] , \\ 
&  [\Tilde{T}_{P_{3, 2}} (E_{\pm , 1} ) \ 
    \Tilde{T}_{P_{3, 2}} (E_{\mp , 2} ) \ 
    \Tilde{T}_{P_{3, 2}} (E_{\pm , 3} )] \\ 
& =[ E_{\pm , 1} \ E_{\mp , 2} \ E_{\pm , 3} ] 
    \left[ 
    \begin{array}{ccc} 
         {\rm cosh}\,t & \mp {\rm sinh}\,t & 0 \\ 
     \mp {\rm sinh}\,t &     {\rm cosh}\,t & 0 \\ 
           0           &       0           & 1 
      \end{array} 
    \right] . 
\end{split} 
\label{P31P32E} 
\end{equation} 
Since 
\begin{equation} 
\begin{array}{lcl} 
\omega_{12} =\dfrac{1}{\sqrt{2}} (E_{+, 1} +E_{-, 1} ), & \ & 
\omega_{34} =\dfrac{1}{\sqrt{2}} (E_{+, 1} -E_{-, 1} ), \\ 
\omega_{13} =\dfrac{1}{\sqrt{2}} (E_{+, 2} +E_{-, 2} ), & \ & 
\omega_{42} =\dfrac{1}{\sqrt{2}} (E_{+, 2} -E_{-, 2} ), \\ 
\omega_{14} =\dfrac{1}{\sqrt{2}} (E_{+, 3} +E_{-, 3} ), & \ & 
\omega_{23} =\dfrac{1}{\sqrt{2}} (E_{+, 3} -E_{-, 3} ), 
  \end{array} 
  \label{omegaE} 
\end{equation} 
applying \eqref{P11P12E} and \eqref{omegaE} to \eqref{Omegasqrt2}, we obtain 
\begin{equation} 
\begin{split} 
      \Tilde{T}_{P_{1,1}} (\Omega ) 
= &  (\cos \phi +1)E_{+, 1} +(\cos \phi -1)E_{-, 1} \\ 
  & - \sin \phi \sin \theta  (E_{+, 2} +E_{-, 2} ) 
    + \sin \phi \cos \theta  (E_{+, 3} -E_{-, 3} ), \\ 
      \Tilde{T}_{P_{1,2}} (\Omega ) 
= &  (\cos \phi +1)E_{+, 1} +(\cos \phi -1)E_{-, 1} \\ 
  & - \sin \phi\,{\rm sinh}\,t(E_{+, 2} -E_{-, 2} ) 
    + \sin \phi\,{\rm cosh}\,t(E_{+, 3} -E_{-, 3} ). 
\end{split} 
\label{P11P12O} 
\end{equation} 
Similarly, applying \eqref{P21P22E}, \eqref{P31P32E} and \eqref{omegaE} 
to \eqref{Omegasqrt2}, we obtain 
\begin{equation} 
\begin{split} 
      \Tilde{T}_{P_{2,1}} (\Omega ) 
= &  (\cos (\phi -\theta )+\cos \theta )E_{+, 1} 
    +(\cos (\phi -\theta )-\cos \theta )E_{-, 1}   \\ 
  & +(\sin (\phi -\theta )-\sin \theta )E_{+, 3} 
    -(\sin (\phi -\theta )+\sin \theta )E_{-, 3} , \\ 
      \Tilde{T}_{P_{2,2}} (\Omega ) 
= &  ((\cos \phi +1)\,{\rm cosh}\,t -\sin \phi\,{\rm sinh}\,t )E_{+, 1} \\ 
  & +((\cos \phi -1)\,{\rm cosh}\,t +\sin \phi\,{\rm sinh}\,t )E_{-, 1} \\ 
  & +((\cos \phi -1)\,{\rm sinh}\,t +\sin \phi\,{\rm cosh}\,t )E_{+, 3} \\ 
  & +((\cos \phi +1)\,{\rm sinh}\,t -\sin \phi\,{\rm cosh}\,t )E_{-, 3} 
\end{split} 
\label{P21P22O} 
\end{equation} 
and 
\begin{equation} 
\begin{split} 
      \Tilde{T}_{P_{3,1}} (\Omega ) 
= &  (\cos \phi +1)\cos \theta E_{+, 1} +(\cos \phi -1)\cos \theta E_{-, 1} \\ 
  & +(\cos \phi +1)\sin \theta E_{+, 2} +(\cos \phi -1)\sin \theta E_{-, 2} \\ 
  & + \sin \phi (E_{+, 3} -E_{-, 3} ), \\ 
      \Tilde{T}_{P_{3,2}} (\Omega ) 
= &  (\cos \phi +1)\,{\rm cosh}\,tE_{+, 1} 
    +(\cos \phi -1)\,{\rm cosh}\,tE_{-, 1} \\ 
  & +(\cos \phi -1)\,{\rm sinh}\,tE_{+, 2} 
    -(\cos \phi +1)\,{\rm sinh}\,tE_{-, 2} \\ 
  & + \sin \phi (E_{+, 3} -E_{-, 3} ). 
\end{split} 
\label{P31P32O} 
\end{equation} 

\section{Tangent vectors of 
\mbox{\boldmath{$SO(3,1)$}}-orbits}\label{sect:tv(3,1)orbits} 

\setcounter{equation}{0} 

By \eqref{P11P12O}, \eqref{P21P22O} and \eqref{P31P32O}, 
we obtain 
\begin{equation} 
\begin{split} 
      \left. \dfrac{\partial}{\partial \theta} \right|_{\theta =0}  
      \Tilde{T}_{P_{1,1}} (\Omega ) 
= &  -\sin \phi (E_{+, 2} +E_{-, 2} ), \\  
      \left. \dfrac{\partial}{\partial t} \right|_{t=0}  
      \Tilde{T}_{P_{1,2}} (\Omega ) 
= &  -\sin \phi (E_{+, 2} -E_{-, 2} ), \\  
      \left. \dfrac{\partial}{\partial \theta} \right|_{\theta =0}  
      \Tilde{T}_{P_{2,1}} (\Omega ) 
= &   \sin \phi (E_{+, 1} +E_{-, 1} ) \\ 
  & -(\cos \phi +1)E_{+, 3} +(\cos \phi -1)E_{-, 3} , \\  
      \left. \dfrac{\partial}{\partial t} \right|_{t=0}  
      \Tilde{T}_{P_{2,2}} (\Omega ) 
= &  -\sin \phi (E_{+, 1} -E_{-, 1} ) \\ 
  & +(\cos \phi -1)E_{+, 3} +(\cos \phi +1)E_{-, 3} , \\ 
      \left. \dfrac{\partial}{\partial \theta} \right|_{\theta =0}  
      \Tilde{T}_{P_{3,1}} (\Omega ) 
= &  (\cos \phi +1)E_{+, 2} +(\cos \phi -1)E_{-, 2} , \\  
      \left. \dfrac{\partial}{\partial t} \right|_{t=0}  
      \Tilde{T}_{P_{3,2}} (\Omega ) 
= &  (\cos \phi -1)E_{+, 2} -(\cos \phi +1)E_{-, 2} . 
\end{split} 
\label{tanvect} 
\end{equation} 
From \eqref{tanvect}, we obtain 

\begin{pro}\label{pro:ts} 
The tangent space $\mathcal{T}_{\Omega} (\mathcal{L} (\Omega ))$ 
of $\mathcal{L} (\Omega )$ at $\Omega$ is spanned by 
\begin{equation*} 
X_{\pm} :=(\sin \phi )E_{\pm , 1} \mp (\cos \phi )E_{\pm , 3} -E_{\mp , 3} , 
   \quad 
Y_{\pm} :=E_{\pm , 2} . 
\label{XpmYpm} 
\end{equation*} 
In addition, 
\begin{itemize} 
\item[{\rm (a)}]{if $\phi =0$ or $\pi$, 
then $\mathcal{T}_{\Omega} (\mathcal{L} (\Omega ))$ is spanned 
by $E_{\pm , 2}$, $E_{\pm , 3}$, i.e., 
$\omega_{13}$, $\omega_{42}$, $\omega_{14}$, $\omega_{23} ;$} 
\item[{\rm (b)}]{if $\phi \not= \pi /2$, 
then the metric $\hat{h}^{\top}$ of $\mathcal{L} (\Omega )$ 
induced by $\hat{h}$ is neutral\/$;$} 
\item[{\rm (c)}]{if $\phi = \pi /2$, 
then $\hat{h}^{\top}$ is degenerate on a two-dimensional subspace $W_0$ 
of $\mathcal{T}_{\Omega} (\mathcal{L} (\Omega ))$ spanned by $X_{\pm}$, 
and therefore every vector of $W_0$ is normal to 
all the vectors of $\mathcal{T}_{\Omega} (\mathcal{L} (\Omega ))$ 
with respect to $\hat{h}^{\top}$.} 
\end{itemize} 
\end{pro} 

\begin{cor}\label{cor:ts} 
If $\phi \not= \pi /2$, then 
\begin{equation*} 
\begin{split} 
X_1 & :=\dfrac{1}{2\cos \phi} 
       (\cos   \phi \sin \phi\,\omega_{12} 
    -(1+\cos^2 \phi )          \omega_{23} 
    -               \sin \phi\,\omega_{34} ), \\ 
X_2 & :=\dfrac{1}{2\cos \phi} 
       (\cos \phi \sin   \phi\,\omega_{12} 
      -2          \cos   \phi\,\omega_{14} 
      +           \sin^2 \phi\,\omega_{23} 
      +           \sin   \phi\,\omega_{34} ), \\ 
Y_1 & :=\omega_{13} , \quad 
Y_2   :=\omega_{42} 
\end{split} 
\end{equation*} 
form a pseudo-orthonormal basis 
of $\mathcal{T}_{\Omega} (\mathcal{L} (\Omega ))$ 
satisfying that $X_1$, $Y_1$ are space-like 
and        that $X_2$, $Y_2$ are time-like. 
\end{cor} 

\section{A special surface in 
an \mbox{\boldmath{$SO(3,1)$}}-orbit}\label{sect:s2idso(3,1)orbits} 

\setcounter{equation}{0} 

Suppose $\phi \not= \pi /2$. 
We set 
\begin{equation} 
\mathcal{S} := \{ \Tilde{T}_{P_{2, 1}}  \circ 
                  \Tilde{T}_{P_{2, 2}} (\Omega ) \ | \ 
                  \theta , t\in \mbox{\boldmath{$R$}} \} . 
\label{surface} 
\end{equation} 
Noticing that $\Tilde{T}_{P_{2, 1}}$ and $\Tilde{T}_{P_{2, 2}}$ are 
commutative, 
we see that $\mathcal{S}$ is a surface in $\mathcal{L} (\Omega )$, i.e., 
a two-dimensional submanifold of $\mathcal{L} (\Omega )$. 
Let $\Tilde{X}_{\pm}$, $\Tilde{Y}_{\pm}$ be vector fields along $\mathcal{S}$ 
given by $\Tilde{X}_{\pm} =\Tilde{T}_{P_{2, 1}} \circ 
                           \Tilde{T}_{P_{2, 2}} (X_{\pm} )$, 
         $\Tilde{Y}_{\pm} =\Tilde{T}_{P_{2, 1}} \circ 
                           \Tilde{T}_{P_{2, 2}} (Y_{\pm} )$. 
Then on $\mathcal{S}$, 
$\Tilde{X}_{\pm}$ (respectively, $\Tilde{Y}_{\pm}$) are 
tangent (respectively, normal) to $\mathcal{S}$. 

Let $\hat{\nabla}$, $\hat{\nabla}^{\top}$ be the Levi-Civita connections 
of $\hat{h}$, $\hat{h}^{\top}$ respectively. 
Since $\Tilde{T}_{P_{2, 1}} (E_{\pm , 2} ) 
      =\Tilde{T}_{P_{2, 2}} (E_{\pm , 2} ) 
      =                      E_{\pm , 2}$ from \eqref{P21P22E}, 
$\Tilde{Y}_{\pm}$ are parallel along $\mathcal{S}$, that is, 
$\hat{\nabla}^{\top}_X \Tilde{Y}_{\pm} =0$ 
for any tangent vector $X$ to $\mathcal{S}$. 
By \eqref{P21P22E}, we obtain 
\begin{equation} 
\begin{split} 
        \Tilde{T}_{P_{2, 1}} (X_{\pm} ) 
= &     \sin (\phi -\theta )E_{\pm , 1} \mp \cos (\phi -\theta )E_{\pm , 3} 
 \pm    \sin        \theta  E_{\mp , 1}  -  \cos \theta         E_{\mp , 3} , 
 \\ 
        \Tilde{T}_{P_{2, 2}} (X_{\pm} ) 
= &    (\sin \phi\,{\rm cosh}\,t -{\rm sinh}\,t)E_{\pm , 1} 
    \mp \cos \phi\,{\rm cosh}\,tE_{\pm , 3} \\ 
  & \mp \cos \phi\,{\rm sinh}\,tE_{\mp , 1} 
      +(\sin \phi\,{\rm sinh}\,t -{\rm cosh}\,t)E_{\mp , 3} . 
\end{split} 
\label{P21P22Xpm}
\end{equation} 
If we set 
\begin{equation*} 
N_{\pm} :=-    \cos \phi E_{\pm , 1} \mp \sin \phi E_{\pm , 3} 
           \pm           E_{\mp , 1} , 
\end{equation*} 
then from \eqref{tanvect} and \eqref{P21P22Xpm}, we obtain 
\begin{equation*} 
\begin{split} 
& \hat{\nabla}_{X_+ +X_-} \Tilde{X}_{\pm} 
 =\left. \dfrac{\partial}{\partial \theta} \right|_{\theta =0} 
  \Tilde{T}_{P_{2, 1}} (X_{\pm} ) 
 =N_{\pm} , \\ 
& \hat{\nabla}_{-X_+ +X_-} \Tilde{X}_{\pm} 
 =\left. \dfrac{\partial}{\partial t} \right|_{t=0} 
  \Tilde{T}_{P_{2, 2}} (X_{\pm} ) 
 =\pm N_{\mp} . 
\end{split} 
\end{equation*} 
We see from Proposition~\ref{pro:ts} 
that $N_{\pm}$ are normal to $\mathcal{L} (\Omega )$ at $\Omega$ 
with respect to $\hat{h}$. 
Since 
\begin{equation*} 
\begin{split} 
   \dfrac{\partial}{\partial \theta} 
   \Tilde{T}_{P_{2, 1}} \circ \Tilde{T}_{P_{2, 2}} (X_{\pm} ) 
& =\Tilde{T}_{P_{2, 1}} \circ \Tilde{T}_{P_{2, 2}} (N_{\pm} ), \\ 
   \dfrac{\partial}{\partial t} 
   \Tilde{T}_{P_{2, 1}} \circ \Tilde{T}_{P_{2, 2}} (X_{\pm} ) 
& =\Tilde{T}_{P_{2, 1}} \circ \Tilde{T}_{P_{2, 2}} (\pm N_{\mp} ), 
\end{split} 
\end{equation*} 
we obtain $\hat{\nabla}^{\top}_X \Tilde{X}_{\pm} =0$ 
for any tangent vector $X$ to $\mathcal{S}$ 
and this means that $\Tilde{X}_{\pm}$ are parallel along $\mathcal{S}$ 
in $\mathcal{L} (\Omega )$. 
Hence we obtain 

\begin{pro}\label{pro:D} 
If $\phi \not= \pi /2$, 
then $\mathcal{S}$ is a time-like surface 
and $\Tilde{X}_{\pm}$, $\Tilde{Y}_{\pm}$ are parallel along $\mathcal{S}$ 
in $\mathcal{L} (\Omega )$. 
\end{pro} 

\begin{rem} 
Suppose $\phi \not= \pi /2$. 
Then by Proposition~\ref{pro:D}, 
$\mathcal{S}$ is flat, and in addition, $\mathcal{S}$ is totally geodesic 
and has flat normal connection in $\mathcal{L} (\Omega )$. 
\end{rem} 

In the case of $\phi = \pi /2$, 
we obtain $N_+ =X_-$ and $N_- =-X_+$. 
Therefore we obtain 

\begin{thm}\label{thm:Dl} 
If $\phi = \pi /2$, 
then there exists a two-dimensional involutive distribution $\mathcal{D}$ 
on $\mathcal{L} (\Omega )$ where $\hat{h}^{\top}$ is degenerate, and 
the covariant derivatives of $\Tilde{X}_{\pm}$ 
$($respectively, $\Tilde{Y}_{\pm} )$ 
along each integral surface $\mathcal{S}$ of $\mathcal{D}$ are 
zero or light-like, and tangent to $\mathcal{S}$ 
$($respectively, vanish$)$. 
\end{thm} 

Noticing $\mathcal{L}_{\pm} 
         =\{ \Tilde{T}_{P} (E_{\pm , 1} ) \ | \ P\in SO(3,1)\}$, 
we will prove 

\begin{thm}\label{thm:gcase} 
If $\phi \not= \pi /2$, 
then $\mathcal{L} (\Omega )$ is homothetic to 
either $\mathcal{L}_+$ or $\mathcal{L}_-$ in $L(\bigwedge^2 E^4_1 )$. 
\end{thm} 

In order to prove Theorem~\ref{thm:gcase}, 
we have only to show the following: 

\begin{pro}\label{pro:forthm} 
Suppose $\phi \not= \pi /2$. 
Let $\mathcal{S}$ be a surface in $\mathcal{L} (\Omega )$ 
given by \eqref{surface}. 
Then $\mathcal{S}$ has a unique element in the form 
of $r(\omega_{12} +\varepsilon \omega_{34} )$ 
for $r\in (0, \sqrt{2} ]$ and $\varepsilon \in \{ +, -\}$. 
\end{pro} 

\vspace{3mm} 

\par\noindent 
\textit{Proof} \ 
If $\phi =0$ or $\pi$, 
then we immediately obtain Proposition~\ref{pro:forthm}. 
In the following, 
suppose $\phi \not= 0$, $\pi /2$, $\pi$. 
The following holds: 
\begin{equation} 
\begin{split} 
       \dfrac{1}{\sqrt{2}} 
       \Tilde{T}_{P_{2, 1}} \circ \Tilde{T}_{P_{2, 2}} (\Omega )
= &  ( \cos (\phi -\theta )\,{\rm cosh}\,t 
      -\sin        \theta  \,{\rm sinh}\,t)\omega_{12} \\ 
  & +(-\sin (\phi -\theta )\,{\rm sinh}\,t 
      +\cos        \theta  \,{\rm cosh}\,t)\omega_{34} \\  
  & +( \cos (\phi -\theta )\,{\rm sinh}\,t 
      -\sin        \theta  \,{\rm cosh}\,t)\omega_{14} \\ 
  & +( \sin (\phi -\theta )\,{\rm cosh}\,t 
      -\cos  \theta        \,{\rm sinh}\,t)\omega_{23} . 
\end{split} 
\label{P21P22Omega} 
\end{equation} 
Both the coefficients of $\omega_{14}$, $\omega_{23}$ 
in the right side of \eqref{P21P22Omega} vanish 
if and only if $\phi$, $\theta$, $t$ satisfy 
\begin{equation} 
\tan 2\theta =\tan \phi , \quad 
\sin \theta\,{\rm cosh}\,t = \cos (\phi -\theta )\,{\rm sinh}\,t. 
\label{ptt} 
\end{equation} 
In addition, the first relation in \eqref{ptt} is 
equivalent to $\theta =\phi /2 +k\pi /2$ 
for an integer $k\in \mbox{\boldmath{$Z$}}$. 
Suppose $\theta =\phi /2 +l\pi$ for $l\in \mbox{\boldmath{$Z$}}$. 
Then the second relation in \eqref{ptt} is rewritten 
into $\tan \phi /2=\,{\rm tanh}\,t$. 
If $\phi \in (0, \pi /2)$, 
then there exists a unique number $t\in \mbox{\boldmath{$R$}}$ 
satisfying $\tan \phi /2=\,{\rm tanh}\,t$; 
if $\phi \in (\pi /2, \pi )$, 
then $\tan \phi /2=\,{\rm tanh}\,t$ does not hold 
for any $t\in \mbox{\boldmath{$R$}}$. 
If $\phi$, $t$ satisfy $\tan \phi /2=\,{\rm tanh}\,t$, 
then from \eqref{P21P22Omega}, we obtain 
\begin{equation} 
 \Tilde{T}_{P_{2, 1}} \circ \Tilde{T}_{P_{2, 2}} (\Omega )
=\sqrt{2} (-1)^l 
 \dfrac{\cos (\phi /2)}{{\rm cosh}\,t} (\omega_{12} +\omega_{34} ). 
 \label{12+34} 
\end{equation} 
Suppose $\theta =\phi /2 +(l+(1/2))\pi$ for $l\in \mbox{\boldmath{$Z$}}$. 
Then the second relation in \eqref{ptt} is rewritten 
into $\cot \phi /2=\,{\rm tanh}\,t$. 
If $\phi \in (0, \pi /2)$, 
then $\cot \phi /2=\,{\rm tanh}\,t$ does not hold 
for any $t\in \mbox{\boldmath{$R$}}$; 
if $\phi \in (\pi /2, \pi )$, 
then there exists a unique number $t\in \mbox{\boldmath{$R$}}$ 
satisfying $\cot \phi /2=\,{\rm tanh}\,t$. 
If $\phi$, $t$ satisfy $\cot \phi /2=\,{\rm tanh}\,t$, 
then from \eqref{P21P22Omega}, we obtain 
\begin{equation} 
 \Tilde{T}_{P_{2, 1}} \circ \Tilde{T}_{P_{2, 2}} (\Omega )
=\sqrt{2} (-1)^l 
 \dfrac{\sin (\phi /2)}{{\rm cosh}\,t} (\omega_{12} -\omega_{34} ). 
 \label{12-34} 
\end{equation} 
From \eqref{12+34} and \eqref{12-34}, 
we immediately obtain Proposition~\ref{pro:forthm}. 
\hfill 
$\square$ 

\vspace{3mm} 

From Proposition~\ref{pro:forthm}, 
we see that there exist $r\in (0, \sqrt{2} ]$ and $\varepsilon \in \{ +, -\}$ 
such that $\Omega_0 :=r(\omega_{12} +\varepsilon \omega_{34} )$ is 
an element of $\mathcal{L} (\Omega )$. 
Therefore $\mathcal{L} (\Omega )$ is an $SO(3, 1)$-orbit through $\Omega_0$. 
Hence we obtain Theorem~\ref{thm:gcase}. 

\section{The stabilizers}\label{sect:stab} 

\setcounter{equation}{0} 

Let $\Omega$ be as in \eqref{Omegasqrt2}. 
Suppose $\phi \not= \pi /2$. 
Then as was already shown in Proposition~\ref{pro:forthm}, 
$\mathcal{L} (\Omega )$ has an element $\Omega_0$ in the form 
of $r(\omega_{12} +\varepsilon \omega_{34} )$. 
The stabilizer of $SO(3, 1)$ at $\Omega_0$ is 
generated by $P_{1, 1}$, $\pm P_{1, 2}$ 
($\theta$, $t\in \mbox{\boldmath{$R$}}$). 
Let $W_+$ (respectively, $W_-$) be a two-dimensional subspace 
of $\mathcal{T}_{\Omega_0} (\mathcal{L} (\Omega ))$ 
generated by $E_{+, 2} +E_{-, 3}$, $E_{-, 2} +E_{+, 3}$ (respectively, 
$E_{+, 2} -E_{-, 3}$, $E_{-, 2} -E_{+, 3}$). 
We will prove 

\begin{pro}\label{pro:invsubsp} 
The proper invariant subspaces 
of $\mathcal{T}_{\Omega_0} (\mathcal{L} (\Omega_0 ))$ by the stabilizer 
are given by $W_{\pm}$. 
\end{pro} 

\vspace{3mm} 

\par\noindent 
\textit{Proof} \ 
Let $\Omega_0$ be an element of $L(\bigwedge^2 E^4_1 )$ 
in the form of $r(\omega_{12} +\varepsilon \omega_{34} )$ 
for $r>0$ and $\varepsilon \in \{ +, -\}$. 
Then $E_{\pm , 2}$, $E_{\pm , 3}$ form a basis of 
the tangent space $\mathcal{T}_{\Omega_0} (\mathcal{L} (\Omega_0 ))$ 
of $\mathcal{L} (\Omega_0 )$ at $\Omega_0$. 
Since the stabilizer of $SO(3, 1)$ at $\Omega_0$ is 
generated by $P_{1, 1}$, $\pm P_{1, 2}$ 
($\theta$, $t\in \mbox{\boldmath{$R$}}$), 
the derivatives $d\Tilde{T}_{P_{1, 1}}$, $d\Tilde{T}_{P_{1, 2}}$ 
of $\Tilde{T}_{P_{1, 1}}$, $\Tilde{T}_{P_{1, 2}}$ respectively 
give linear transformations 
of $\mathcal{T}_{\Omega_0} (\mathcal{L} (\Omega_0 ))$. 
Let $W$ be an invariant subspace 
of $\mathcal{T}_{\Omega_0} (\mathcal{L} (\Omega_0 ))$ by the stabilizer. 
Let 
$$E=aE_{+, 2} +bE_{-, 2} +cE_{+, 3} +dE_{-, 3}$$ 
be an element of $W$. 
Then we have 
\begin{equation*} 
 d\Tilde{T}_{P_{1, 1}} \circ d\Tilde{T}_{P_{1, 2}} (E) 
=(E_{+, 2} \ E_{-, 2} \ E_{+, 3} \ E_{-, 3} )\mbox{\boldmath{$a$}} , 
\end{equation*} 
where 
\begin{equation*} 
  \mbox{\boldmath{$a$}} 
:=\cos \theta\,{\rm cosh}\,t\,\mbox{\boldmath{$a$}}_1 
 +\cos \theta\,{\rm sinh}\,t\,\mbox{\boldmath{$a$}}_2 
 +\sin \theta\,{\rm cosh}\,t\,\mbox{\boldmath{$a$}}_3 
 +\sin \theta\,{\rm sinh}\,t\,\mbox{\boldmath{$a$}}_4 
\end{equation*} 
and 
\begin{equation*} 
  \mbox{\boldmath{$a$}}_1 
:=\left[ \begin{array}{c} 
          a \\ 
          b \\ 
          c \\ 
          d 
           \end{array} 
  \right] , \ 
  \mbox{\boldmath{$a$}}_2 
:=\left[ \begin{array}{c} 
          d \\ 
          c \\ 
          b \\ 
          a 
           \end{array} 
  \right] , \ 
  \mbox{\boldmath{$a$}}_3 
:=\left[ \begin{array}{c} 
         -c \\ 
          d \\ 
          a \\ 
         -b 
           \end{array} 
  \right] , \ 
  \mbox{\boldmath{$a$}}_4 
:=\left[ \begin{array}{c} 
         -b \\ 
          a \\ 
          d \\ 
         -c 
           \end{array} 
  \right] . 
\end{equation*} 
We set $A:=[\mbox{\boldmath{$a$}}_1 \ 
            \mbox{\boldmath{$a$}}_2 \ 
            \mbox{\boldmath{$a$}}_3 \ 
            \mbox{\boldmath{$a$}}_4 ]$. 
Then we have 
\begin{equation*} 
\det A=-((a-d)^2 +(b-c)^2 )((a+d)^2 +(b+c)^2 ). 
\end{equation*} 
Therefore we see that if $\det A\not= 0$, 
then $W$ coincides with $\mathcal{T}_{\Omega_0} (\mathcal{L} (\Omega_0 ))$. 
Suppose $\det A=0$. 
Then we have either $(a, b)=(d, c)$ or $(a, b)=(-d, -c)$. 
Therefore we see that if $W$ has a nonzero element, 
then $W=W_+$ or $W_-$, where $W_{\pm}$ are as in the beginning  
of this section. 
Hence we obtain Proposition~\ref{pro:invsubsp}. 
\hfill 
$\square$ 

\vspace{3mm} 

\begin{rem} 
As in \cite{ando}, \cite{ando'}, 
there exist an almost complex structure $\mathcal{I}$ and 
an almost paracomplex structure $\mathcal{J}$ on $\mathcal{L} (\Omega_0 )$ 
parallel with respect to the Levi-Civita connection $\hat{\nabla}^{\top}$. 
We see that $W_{\pm}$ are invariant by $\mathcal{I}$, $\mathcal{J}$. 
\end{rem} 

In the following, suppose $\phi = \pi /2$. 
Then we have 
\begin{equation} 
\Omega =E_{+, 1} -E_{-, 1} +E_{+, 3} -E_{-, 3} . 
\label{Omegaphipi/2} 
\end{equation} 
For $k=1, 2, 3$, 
we denote $P_{k, 1}$, $P_{k, 2}$ 
by $P_{k, 1, \theta}$, $P_{k, 2, t}$ respectively. 
For $t\in \mbox{\boldmath{$R$}}$, we set 
\begin{equation} 
\theta (t):= \sin^{-1} ({\rm tanh}\,t), \quad 
s      (t):=-\log      ({\rm cosh}\,t) 
\label{ths} 
\end{equation} 
and 
\begin{equation} 
U_t :=P_{3, 2, t} P_{1, 1, -\theta (t)} P_{2, 2, s(t)} , \quad 
V_t :=P_{1, 2, t} P_{3, 1,  \theta (t)} P_{2, 2, s(t)} . 
\label{utvt} 
\end{equation} 
We will prove 

\begin{thm}\label{thm:stab} 
The stabilizer of $SO(3, 1)$ at $\Omega$ is generated by 
$U_t$, $V_t$ in \eqref{utvt} $(t\in \mbox{\boldmath{$R$}})$ 
and each element of the stabilizer fixes each vector 
of the two-dimensional subspace $W_0$ 
of $\mathcal{T}_{\Omega} (\mathcal{L} (\Omega ))$ spanned by $X_{\pm}$. 
In addition, 
\begin{itemize} 
\item[{\rm (a)}]{a one-dimensional subspace 
of $\mathcal{T}_{\Omega} (\mathcal{L} (\Omega ))$ is invariant 
by the stabilizer if and only if it is contained in $W_0 ;$} 
\item[{\rm (b)}]{$W_0$ is a unique invariant two-dimensional subspace 
of $\mathcal{T}_{\Omega} (\mathcal{L} (\Omega ))$ 
by the stabilizer\/$;$} 
\item[{\rm (c)}]{a three-dimensional subspace 
of $\mathcal{T}_{\Omega} (\mathcal{L} (\Omega ))$ is invariant 
by the stabilizer if and only if it contains $W_0$.} 
\end{itemize} 
\end{thm} 

\vspace{3mm} 

\par\noindent 
\textit{Proof} \ 
From \eqref{P21P22O} with $\phi =\pi /2$ and \eqref{Omegaphipi/2}, 
we have $\Tilde{T}_{P_{2, 2, s}} (\Omega )=e^{-s} \Omega$ 
for any $s\in \mbox{\boldmath{$R$}}$. 
Therefore by \eqref{P11P12O} with $\phi =\pi /2$, we have 
\begin{equation} 
\begin{split} 
& \Tilde{T}_{P_{1, 1, \theta}} \circ \Tilde{T}_{P_{2, 2, s}} (\Omega ) \\ 
& =e^{-s} (E_{+, 1} -E_{-, 1} -\sin \theta (E_{+, 2} +E_{-, 2} ) 
                              +\cos \theta (E_{+, 3} -E_{-, 3} )) 
\end{split} 
\label{11th22s} 
\end{equation} 
for any $\theta \in \mbox{\boldmath{$R$}}$. 
By \eqref{P31P32E} and \eqref{11th22s}, we have 
\begin{equation} 
\begin{split} 
& \Tilde{T}_{P_{3, 2, t}} \circ \Tilde{T}_{P_{1, 1, \theta}} 
                          \circ \Tilde{T}_{P_{2, 2, s}} (\Omega ) \\ 
& =e^{-s} (({\rm cosh}\,t +\sin \theta\,{\rm sinh}\,t)(E_{+, 1} -E_{-, 1} ) \\ 
&            \qquad \quad 
          -({\rm sinh}\,t +\sin \theta\,{\rm cosh}\,t)(E_{+, 2} +E_{-, 2} ) 
          +                \cos \theta                (E_{+, 3} -E_{-, 3} )) 
\end{split} 
\label{32t11th22s} 
\end{equation} 
for any $t\in \mbox{\boldmath{$R$}}$. 
Let $\theta (t)$, $s(t)$ be as in \eqref{ths}. 
Then noticing 
\begin{equation*} 
 e^{s(t)} 
={\rm cosh}\,t +\sin (-\theta (t))\,{\rm sinh}\,t 
=\cos \theta (t) 
=\dfrac{1}{{\rm cosh}\,t} 
\end{equation*} 
and 
\begin{equation*} 
{\rm sinh}\,t +\sin (-\theta (t))\,{\rm cosh}\,t=0, 
\end{equation*} 
we obtain 
\begin{equation*} 
 \Tilde{T}_{U_t} (\Omega ) 
=\Tilde{T}_{P_{3, 2, t}} \circ \Tilde{T}_{P_{1, 1, -\theta (t)}} 
                         \circ \Tilde{T}_{P_{2, 2, s(t)}} (\Omega ) 
=\Omega . 
\end{equation*} 
Similarly, we obtain 
\begin{equation*} 
 \Tilde{T}_{V_t} (\Omega ) 
=\Tilde{T}_{P_{1, 2, t}} \circ \Tilde{T}_{P_{3, 1, \theta (t)}} 
                         \circ \Tilde{T}_{P_{2, 2, s(t)}} (\Omega ) 
=\Omega . 
\end{equation*} 
Therefore $U_t$, $V_t$ ($t\in \mbox{\boldmath{$R$}}$) are elements 
of the stabilizer of $SO(3, 1)$ at $\Omega$. 
We set 
\begin{equation*} 
U'_x :=\left[ \begin{array}{cccc} 
               1 &         x       & 0 &         x       \\ 
              -x & 1-\frac{x^2}{2} & 0 &  -\frac{x^2}{2} \\ 
               0 &         0       & 1 &         0       \\ 
               x &   \frac{x^2}{2} & 0 & 1+\frac{x^2}{2}  
                \end{array} 
        \right] , \quad 
V'_x :=\left[ \begin{array}{cccc} 
               1 &         0       &  0 &         0       \\ 
               0 & 1-\frac{x^2}{2} & -x &  -\frac{x^2}{2} \\ 
               0 &         x       &  1 &         x       \\ 
               0 &   \frac{x^2}{2} &  x & 1+\frac{x^2}{2}  
                \end{array} 
        \right] . 
\end{equation*} 
Then we have $U'_{{\rm tanh}\,t} =U_t$, $V'_{{\rm tanh}\,t} =V_t$. 
We see that $\{ U'_x \ | \ x\in \mbox{\boldmath{$R$}} \}$, 
            $\{ V'_x \ | \ x\in \mbox{\boldmath{$R$}} \}$ are 
one-parameter subgroups of $SO(3, 1)$ 
which are different from each other. 
Since the dimension of the stabilizer of $SO(3, 1)$ at $\Omega$ is equal to 
two, it is generated by 
$U_t$, $V_t$ in \eqref{utvt} $(t\in \mbox{\boldmath{$R$}})$. 
We set 
\begin{equation*} 
\Omega':=E_{+, 1} +E_{-, 1} -E_{+, 3} -E_{-, 3} . 
\end{equation*} 
Then $U_t$, $V_t$ fix $\Omega'$. 
Since $\Omega$, $\Omega'$ form a basis of $W_0$, 
each element of the stabilizer fixes each vector of $W_0$. 
Noticing 
\begin{equation*} 
\begin{split} 
  \Tilde{T}_{U_t} (E_{+, 2} ) 
& =E_{+, 2} -\dfrac{{\rm tanh}\,t}{2} (\Omega -\Omega' ), \\ 
  \Tilde{T}_{U_t} (E_{-, 2} ) 
& =E_{-, 2} -\dfrac{{\rm tanh}\,t}{2} (\Omega +\Omega' ) 
\end{split} 
\end{equation*} 
and 
\begin{equation*} 
\begin{split} 
  \Tilde{T}_{V_t} (E_{+, 2} ) 
& =E_{+, 2} -\dfrac{{\rm tanh}\,t}{2} (\Omega +\Omega' ), \\ 
  \Tilde{T}_{V_t} (E_{-, 2} ) 
& =E_{-, 2} +\dfrac{{\rm tanh}\,t}{2} (\Omega -\Omega' ), 
\end{split} 
\end{equation*} 
we obtain (a), (b), (c) of Theorem~\ref{thm:stab}. 
\hfill 
$\square$ 

\begin{rem}
By direct computations, 
we see that $U'_x$ and $V'_y$ are commutative 
for arbitrary $x$, $y\in \mbox{\boldmath{$R$}}$. 
\end{rem} 

\section{The \mbox{\boldmath{$r$}}-slice 
of \mbox{\boldmath{$L(\bigwedge^2 E^4_1 )$}}}\label{sect:rslice} 

\setcounter{equation}{0} 

For an equivalence class $[e]$ in $\mathcal{B} (E^4_1 )$ and each $r>0$, 
the $r$-\textit{slice\/} of $L(\bigwedge^2 E^4_1 )$ given by 
$$L_r (\textstyle\bigwedge^2 E^4_1 , [e]) 
:=\left\{ \left. 
  \Omega \in L(\textstyle\bigwedge^2 E^4_1 ) 
  \ \right| \ 
A(\Omega ) =B(\Omega ) =r^2 
  \right\}$$ 
is well-defined. 
The $1$-slice of $L(\bigwedge^2 E^4_1 )$ coincides 
with $L_1 (\bigwedge^2 E^4_1 , [e])$ given in Section~\ref{sect:lc}. 
The $r$-slice of $L(\bigwedge^2 E^4_1 )$ depends on 
the choice of an equivalence class $[e]$ in $\mathcal{B} (E^4_1 )$. 
In the following, we fix an equivalence class $[e]$ 
and we simply denote the $r$-slice of $L(\bigwedge^2 E^4_1 )$ 
by $L_r (\bigwedge^2 E^4_1 )$. 
Let $\Omega$ be as in \eqref{Omegasqrt2} 
and $\mathcal{L} (\Omega )$ as in \eqref{HOmega}. 
We set $\mathcal{L}_r (\Omega ) 
      :=\mathcal{L}   (\Omega )\cap L_r (\bigwedge^2 E^4_1 )$. 
For $\phi \in [0, \pi ]$, we set 
$$r(\phi ) 
  :=\inf \{ r>0 \ | \ \mathcal{L}_r (\Omega )\not= \emptyset \} .$$ 
For $\phi \in [0, \pi ]\setminus \{ \pi /2\}$, we set 
\begin{equation*} 
t(\phi )
:=\left\{ 
  \begin{array}{cl} 
   {\rm \tanh}^{-1} (\tan (\phi /2)) & (\phi \in [0, \pi /2)), \\ 
   {\rm \tanh}^{-1} (\cot (\phi /2)) & (\phi \in ( \pi /2 , \pi ]). 
    \end{array} 
  \right. 
\end{equation*} 
We have $t(\pi -\phi )=t(\phi )$ for $\phi \in [0, \pi /2)$. 
We will prove 

\begin{thm}\label{thm:slice1} 
For $\phi \in [0, \pi ]\setminus \{ \pi /2\}$, 
$r(\phi )$ is positive and given by 
\begin{equation} 
r(\phi ) 
 =\left\{ 
  \begin{array}{cl} 
   \sqrt{2}\,\dfrac{\cos (\phi /2)}{{\rm cosh}\,t(\phi )} & 
  (\phi \in [0, \pi /2)), \\ 
   \sqrt{2}\,\dfrac{\sin (\phi /2)}{{\rm cosh}\,t(\phi )} & 
  (\phi \in (\pi /2 , \pi ]). 
    \end{array} 
  \right. 
\label{rphi} 
\end{equation} 
In particular, 
a function $r$ of $\phi$ satisfies $r(\pi -\phi )=r(\phi )$ 
and $\displaystyle\lim_{\phi \rightarrow \pi /2} r(\phi )=0$, 
and it is monotonically decreasing on $[0, \pi /2)$. 
In addition, 
\begin{itemize} 
\item[{\rm (a)}]{$\mathcal{L}_r (\Omega )$ is not empty, and diffeomorphic 
to $\mbox{\boldmath{$R$}}P^3$ for any $r>r(\phi );$} 
\item[{\rm (b)}]{$\mathcal{L}_{r(\phi )} (\Omega )$ is not empty, 
and diffeomorphic to $S^2$.} 
\end{itemize} 
\end{thm} 

\vspace{3mm} 

\par\noindent 
\textit{Proof} \ 
We denote by $r_0$ the number given by the right side of \eqref{rphi}. 
Then noticing \eqref{12+34} and \eqref{12-34}, 
we see that $\mathcal{L}_{r_0} (\Omega )$ is not empty, 
diffeomorphic to $S^2$ and given by 
\begin{equation*} 
 \mathcal{L}_{r_0} (\Omega ) 
=\{ r_0 ((\omega_{23} \ \omega_{31} \ \omega_{12} ) 
        + \varepsilon 
         (\omega_{14} \ \omega_{24} \ \omega_{34} ))\mbox{\boldmath{$a$}} 
          \ | \ 
          \langle \mbox{\boldmath{$a$}} , \mbox{\boldmath{$a$}} \rangle =1 \} 
\end{equation*} 
for $\varepsilon \in \{ +, -\}$. 
If $r>r_0$, 
then we obtain $\mathcal{L}_r (\Omega )\not= \emptyset$ 
using $\Tilde{T}_{P_{1, 2}}$, $\Tilde{T}_{P_{2, 2}}$, $\Tilde{T}_{P_{3, 2}}$, 
and noticing that $\mathcal{L}_r (\Omega )$ is represented as 
\begin{equation*} 
 \mathcal{L}_r (\Omega )
=\{ r((\omega_{23} \ \omega_{31} \ \omega_{12} )U\mbox{\boldmath{$a$}} 
     +(\omega_{14} \ \omega_{24} \ \omega_{34} )U\mbox{\boldmath{$b$}} )
       \ | \ U\in SO(3)\} 
\end{equation*} 
with $\mbox{\boldmath{$b$}} \not= \pm \mbox{\boldmath{$a$}}$, 
we see that $\mathcal{L}_r (\Omega )$ is diffeomorphic 
to $\mbox{\boldmath{$R$}}P^3$. 
On a neighborhood of an element of $\mathcal{L}_{r_0} (\Omega )$ 
in $\mathcal{L} (\Omega )$, 
there exist no elements of $\mathcal{L}_r (\Omega )$ with $r<r_0$. 
Therefore we obtain $\mathcal{L}_r (\Omega )=\emptyset$ for $r\in (0, r_0 )$. 
Hence we obtain Theorem~\ref{thm:slice1}. 
\hfill 
$\square$ 

\vspace{3mm} 

For $r>0$ and $\varepsilon \in \{ +, -\}$, 
an element $r(\omega_{12} +\varepsilon \omega_{34} )$ is 
contained in a neutral $SO(3, 1)$-orbit. 
This means that an $SO(3, 1)$-orbit 
with a two-dimensional involutive distribution 
where $\hat{h}^{\top}$ is degenerate does not 
contain an element in the form 
of $r(\omega_{12} +\varepsilon \omega_{34} )$ 
for $r>0$ and $\varepsilon \in \{ +, -\}$. 
By \eqref{P21P22O} with $\phi =\pi /2$, we obtain $r(\pi /2)=0$ 
and we see that $\Omega$ with $\phi =\pi /2$ satisfies 
$\mathcal{L}_r (\Omega )\not= \emptyset$ for any $r>0$. 
Hence we obtain 

\begin{thm}\label{thm:slice2} 
Suppose $\phi =\pi /2$. 
Then $r(\phi )=0$, and $\mathcal{L}_r (\Omega )$ is not empty, 
and diffeomorphic to $\mbox{\boldmath{$R$}}P^3$ for any $r>0$. 
\end{thm} 

\begin{rem} 
Let $\Omega$ be as in \eqref{Omegasqrt2}. 
Then $\phi \in [0, \pi ]$ depends on the choice of 
an equivalence class $[e]$ in $\mathcal{B} (E^4_1 )$. 
However, noticing Theorem~\ref{thm:slice1} and Theorem~\ref{thm:slice2}, 
we see that whether $\phi$ is equal to $\pi /2$ or not does not depend 
on the choice of $[e]$. 
\end{rem} 

\section*{Acknowledgments} 

The author is grateful to Professors Osamu Ikawa and Takahiro Hashinaga 
for valuable discussions and comments. 
This work was supported by JSPS KAKENHI Grant Number JP21K03228.

\vspace{4mm} 

\par\noindent 
\footnotesize{Faculty of Advanced Science and Technology, 
              Kumamoto University \\ 
              2--39--1 Kurokami, Chuo-ku, Kumamoto 860--8555 Japan} 

\par\noindent  
\footnotesize{E-mail address: andonaoya@kumamoto-u.ac.jp} 

\end{document}